\documentclass[11pt]{amsart}

\usepackage{amscd,amsmath, amssymb, fancyhdr}

\newcommand{\version}{version 1.01,\ \   January 06, 2010}

\renewcommand{\leftmark}{{\scriptsize  LCK metrics on elliptic principal bundles}}

\newtheorem{thm}{Theorem}

\newtheorem{lem}{Lemma}

\newtheorem{coro}{Corollary}

\def\C{\mathbb C}

\def\P{\mathbb P}

\def\R{\mathbb R}

\def\Z{\mathbb Z}

\def\defq{\stackrel{def}{=}}

\def\ra{\rightarrow }

\def\noi{\noindent}

\title{LCK metrics on elliptic principal bundles}

\author{Victor Vuletescu}

\date{}

\begin{document}

\begin{abstract}
For elliptic principal bundles $\pi:X\ra B$ over K\"ahler manifolds it was 
shown by Blanchard that $X$ has a K\"ahler metric if and only both Chern classes (with real coefficients) of $\pi$ vanish. For some elliptic principal bundles, when  the span of these Chern classes is 1-dimensional, it was shown by Vaisman that $X$ carry locally conformally K\"ahler (LCK, for short)  metrics. We show that in the case when the Chern classes are linearly independent, $X$ carries no LCK metric. 
\end{abstract}

\maketitle

\section{Introduction.}

\subsection{Motivation.}

Elliptic principal bundles, or, more generally, principal torus bundles  (we will recall imediately the definitions) were always an excellent ``reservoir" for interesting examples in complex geometry. A short list would certainly include: 

- the very first non-K\"ahler manifolds (H. Hopf, 1948) - some Hopf surfaces - which are such  bundles over $\P^1$; 

- the first examples of non-K\"ahler simply connected manifolds (Calabi - Eckmann manifolds) are elliptic bundles over a product of projective spaces $\P^r\times \P^s$, 

- the first examples of manifolds with non-closed global homolorphic forms (Iwasawa manifolds) are again elliptic principal bundles over  abelian surfaces,

- examples of manifolds for which the Fr\"ohlicher spectral sequence degenerates arbitrarly high (R\"ollenske, \cite{Ro}).

 On the other hand, these kind of manifolds are abundant; for instance, in complex dimension 2 such compact surfaces exists in each class in Kodaira classification (except for class $VII_{>0}$).

 A natural question one may ask about this type of manifolds is about the kind of hermitian metrics they can carry. A classical result of Blanchard (cf. \cite{Bla}) states that if $\pi:X\ra B$ is an elliptic principal bundle over a compact K\"ahler manifold, then $X$ carries a K\"ahler metric if and only if both Chern classes (with real coefficients) of the bundle vanish. On the other hand, on a rather large class of elliptic principal bundles which are non-K\"ahler, one can show the existence of a locally conformally K\"ahler metric.




\subsection{Chern classes of elliptic principal bundles.}

In this section, we quickly recall some basic facts about elliptic principal bundles: more details can be found for instance in \cite{Bas} or in \cite{Ho}.
Let $B$ be a differentiable manifold and $E$ a 2-dimensional real torus. The set of isomorphism classes of principal  bundles over $B$ with fiber $E$ is classified by the cohomology group
 $H^1(B, {\mathcal C}_B(E))$ where 
${\mathcal C}_B(E)$ is a ({\em ad-hoc}) notation for the sheaf of germs of differentiable, $E-$valued functions defined on (open subsets of) $B.$
 As $E\simeq S^1\times S^1$, letting $\Z_B$ (respectively ${\mathcal C}_B(\C)$) for the sheaf of germs of differentiable integer-valued (respectively complex-valued) functions, from the exact sequence 

$$
0
\ra
\Z^{\oplus 2}_B
\ra
{\mathcal C}_B(\C)
\ra
{\mathcal C}_B(E)
\ra
0
$$
and the  fact that ${\mathcal C}_B(\C)$ is a fine sheaf, we see
$$H^1(B,{\mathcal C}_B(E))\simeq H^2(B, \Z_B)^{\oplus 2}.$$

\noi The image of a given $E-$principal bundle $X\ra B$ over $B$ under the above isomorphism will be denoted (keeping the notations in \cite{Bas}) by
$(c_1'(X),c_1''(X));$ its components will be called {\em the Chern classes of $X$}.

An alternative way of defining the Chern classes is as follows. 
As translations of $E$ act trivially in the cohomology of $E$, we see that for any  principal bundle 
$\pi:X\ra B$ with fiber $E$ one  has ${\mathcal R}^i\pi_*(\Z_X)\simeq {\mathbb Z}_B\otimes H^i(E, \Z)$ for all $i.$
Fix $\alpha, \beta$ generators of $H^2(E, \Z)$ and consider the spectral sequence

$$
E_2^{pq}=H^q(B, {\mathcal R}^p\pi_*(\Z_X))\Rightarrow H^{p+q}(B, \Z_B).
$$
Then the images of $\alpha, \beta$ 
 under the differential

$$d:H^0(B, {\mathcal R}^1\pi_*(\Z_X))\ra
 H^2(B, \pi_*(\Z_X))$$
are the above Chern classes (modulo a possible twist by an automorphism of $H^1(E, \Z),$ that is, modulo action by an element in $SL_2(\Z)$).

In the next lemma we gather some basic facts needed further; most likely, its content is well-known, but as we were not able to find a precise reference, we include a proof below. Notice that Chern classes here are viewed in $H^2(B, \R)$ under the natural map $H^2(B, \Z)\ra H^2(B, \R).$

\begin{lem}\label{lema2}
Let $\pi:X\ra B$ be a principal bundle with fiber a 2-dimensional torus $E.$

a) If at least one of the Chern classes of it is nontrivial in $H^2(B, \R)$ then  the homology class (with $\R$-coefficients)  of a fiber vanishes.

b) If the two Chern classes are linealy independent in $H^2(B, \R)$ then the natural map
$\pi^*:H^1(B, \R)\ra H^1(X, \R)$
is an isomorphism.

\end{lem}

\begin{proof}
 a) Let ${\mathbb R}_X$ (resp. ${\mathbb R}_B$) denote the sheaf of locally constant
 functions on $X$ (resp. $B$). Consider the spectral sequence

$$
E_2^{pq}=H^q(B, {\mathcal R}^p\pi_*(\R_X))\Rightarrow H^{p+q}(B, \R_B).
$$
\noi This induces an exact sequence

$$
H^2(X, \R)\stackrel{i^*}{\ra}H^0(B, {\mathcal R}^2\pi_*(\R_X))\stackrel{d}{\ra}H^2(B, {\mathcal R}^1\pi_*(\R_X))$$

\noi which, under canonical identifications, becomes

\begin{eqnarray}\label{exc2}
H^2(X, \R)
\stackrel{i^*}{\ra}
H^2(E,\R) \stackrel{d}{\ra}
H^2(B,\R)\otimes H^1(E, \R).
\end{eqnarray}
We must show that the natural map $H_2(E, \R)\ra
 H_2(X, \R)$ induced by the inclusion of a fiber is zero.
This map is the dual of the map $i^*$ form (\ref{exc2}).
But $H^2(E,\R)$ is generated by $\alpha\wedge\beta$ and as differentials in the above spectral sequence are  multiplicative, we see

$$
d(\alpha\wedge\beta)=c_1'(X)\otimes \beta+c_1''(X)\otimes \alpha$$
Now if at least one of $c_1'(X), c_1''(X)$ is non-vanishing, we see the map $d$ is also nonvanishing; since $H^2(E, \R)$ is one dimensional, we get $d$ is actually injective, so $i^*$ is the null map.

b). The proof is similar, only at this time we look at another exact sequence induced by the spectral sequence above, namely (again under canonical identifications) at
$$
0
\ra
H^1(B, \R)
\stackrel{\pi^*}{\ra}
H^1(X, \R)
\stackrel{\delta}{\ra}
H^1(E, \R)
\stackrel{c}{\ra}
H^2(B, \R)$$
The map $c$ acts by $c(\alpha)=c_1', c(\beta)=c_1''$ so our hypothesis implies $c$ is injective, hence $\delta$ is the null map thus $\pi^*$ is surjective.

\end{proof}

\medskip

If $X, B$ are complex manifolds and $E$ is an elliptic curve (viewed as a complex 1-dimensional Lie group), a holomorphic principal bundle map  $\pi:X\ra B$ with fiber $E$ will be called {\em elliptic principal bundle} for short. 
From  a) of the above Lemma and since no proper compact complex submanifold of a K\"ahler manifold can be homologous to zero, we have:

\begin{coro}\label{corob} If $E$ is an elliptic curve and $\pi:X\ra B$ is an elliptic principal bundle such that at least one of the Chern classes is non-vanishing (with real coefficients) then $X$ carries no K\"ahler metric.

\end{coro}

In particular, we rediscover one implication in Blanchard's theorem - but with no hypothesis on the base $B.$

\subsection{Locally conformally K\"ahler structures.}
We next recall some basic facts about locally conformally K\"ahler structures; more details can be found in the monograph \cite{DrOr}.
By definition, given a complex manifold $X$, a hermitian metric $g$ on $X$ is called {\it locally conformally K\"ahler} (LCK for short) if $X$ can be covered by open subsets $X=\displaystyle\bigcup_{i\in I}U_i$ with the property that for each $i\in I$ there exists a K\"ahler metric $g_i$ defined on $U_i$ such that $g_{\vert U_i}$ is conformally equivalent to $g_i,$ {\em i.e.} there exists some real-valued smooth function $f_i$ defined on $U_i$ such that $g_{\vert U_i}=e^{f_i}g_i.$ The definition is equivalent to the following. Let $\Omega$ be the K\"ahler form of $g;$  then $g$ is LCK iff there exists a closed $1-$form $\theta$ such that 
$$d\Omega=\theta\wedge \Omega.$$
Moreover, $g$ is globally conformal to a K\"ahler metric iff $\theta$ is exact.

Historically, the first example of an LCK metric was given by Vaisman in 1976 in \cite{Va}; we include it here very briefly. Let $n\geq 1$ and let $W\defq\C^n\setminus\{0\}.$ Let $\Z$ act on $W$ by 
$$l\circ(z_1,\dots, z_n)=(2^lz_1,\dots,2^lz_n),l\in \Z.$$
Then the action is fixed-point-free and properly discontinous, and the quotient $X\defq W/\Z$ is  compact. One immediately checks that the hermitian metric
$$g=\frac{1}{\vert z_1\vert^2+\dots+\vert z_n\vert ^2}(dz_1d{\overline z}_1+\dots +dz_nd{\overline z}_n)$$
on $W$ is $\Z-$invariant and defines an LCK metric on $X.$

On the other hand, $X$ has an obvious holomorphic projection
$\pi:X\ra {\mathbb P}^{n-1}(\C);$ one immediately checks that $\pi$ is actually an elliptic principal bundle, whose Chern classes are $(\alpha, 0)$ where $\alpha$ is the Poincar\'e dual of a  hiperplane section in ${\mathbb P}^{n-1}(\C).$

As a consequence, we get the following 
\begin{coro}
Let $B$ be a complex projective manifold. Then there exists an elliptic principal bundle $X\ra B$ such that $X$ has no K\"ahler metric but carries an LCK metric. For this bundle one of the Chern classes vanishes. 
\end{coro}
Indeed, just embed $B$ into some projective space and take the restriction of the map $\pi$ above to it.


\section{The main result.}

We begin by recalling, for the sake of simplicity, the following well-known statement:
\begin{lem}\label{lema1}

Let $(X, \Omega)$ be an LCK manifold, $\theta$ the associated Lee form of $\Omega$ and let $\omega$ be a $1-$form on $X$ which is cohomologous to $\theta$. Then there exists a metric $\Omega'$ on $X$ which is conformally equivalent to $\Omega$ (in particular $\Omega'$ is LCK too) such that the Lee form of $\Omega'$ is $\omega.$ 
\end{lem}

\noi We are now ready to prove the main result of this note.

\begin{thm}

Let $X,B$ be  compact complex manifolds, $X\stackrel{\pi}{\ra} B$ an elliptic principal bundle with fiber $E$.  If  the Chern classes of the bundle $X\stackrel{\pi}{\ra}B$ are linearly independent in $H^2(B, \R)$ then $X$ carries no locally conformally K\"ahler structure.

\end{thm}

\begin{proof}
Assume $X$ carries an LCK metric $\Omega$ with  Lee form$\theta.$ As the Chern classes are independent, we get that $\pi^*:H^1(B, \R)\ra H^1(X, \R)$ is an isomorphism.
Hence there exists a $1-$form $\eta$ on $B$ such that $\theta$ and $\pi^*(\eta)$ are cohomologous. By Lemma \ref{lema1} we may assume that $\theta=\pi^*(\eta).$ This yields 

\begin{eqnarray}\label{fun}d\Omega=\pi^*(\eta)\wedge \Omega
.\end{eqnarray}

Now for each $b\in B$ let $\lambda(b)\defq \int_{E_b}\Omega=vol_{\Omega}(E_b);$ clearly $\lambda$ is a differentiable map $\lambda:B\ra \R_{>0}.$

Hence, replacing  $\Omega$ with $\frac{1}{\lambda}\Omega$ we can assume that relation (\ref{fun}) holds good and also that the volume of the fibers with respect to $\Omega$ is constant.

Let now $[\gamma]\in \pi_1(B)$ be arbitrary and let $\gamma\in [\gamma]$ be a smooth representative. Let $M\defq\pi^{-1}(\gamma);$ we see $M$ is a smooth real $3-$submanifold of $X.$ 
We have
$$
0=\int_Md\Omega=\int_M\pi^*(\eta)\wedge\Omega.$$
But from the independence of $\int_{E_b}\Omega$ on $b\in B$ and from Fubini's theorem, we get furthermore

$$
0
=\left(\int_\gamma\eta\right)\left(\int_{E_b}\Omega\right)
.$$
We obtain $\int_\gamma\eta=0$ and hence $\eta$ is exact, as $[\gamma]$ was arbitrary. We derive that $\Omega$ is globally conformally K\"ahler, so $X$ would admit a K\"ahler metric; but this is impossible by Corollary \ref{corob}.
\end{proof}


{\small

\noindent {\sc Victor Vuletescu\\
University of Bucharest, Faculty of Mathematics and Informatics, \\14 
Academiei str., 70109 Bucharest, Romania. }}
\end{document}